\documentclass[12pt,intlim,righttag]{article}
\usepackage{amsmath,amsthm}
\usepackage{amssymb}
\usepackage{amsfonts}

\newcommand{\la}{\langle}
\newcommand{\ra}{\rangle}

\newcommand\beq{\begin{equation}}
\newcommand\eeq{\end{equation}}

\theoremstyle{Theorem}

\theoremstyle{corollary}

\theoremstyle{remark}

\theoremstyle{definition}

\begin{document}
\title{Necessary and sufficient conditions for the uniform integrability of the stochastic exponential}

\author{B. Chikvinidze }

\date{~}
\maketitle

\begin{center}

  Institute of Cybernetics of Georgian Technical University,
  
  \
  \\       
  Georgian-American University, Business school, 8 M. Aleksidze Srt., \\
          Tbilisi 0160, Georgia 

\          
\\
   E-mail: beso.chiqvinidze@gmail.com \\

\end{center}

\begin{abstract}


{\bf Abstract}
\;We establish necessary and sufficient conditions for the uniform integrability of the stochastic exponential $\cal E(M)$.

\bigskip

\noindent {\it  2000 Mathematics Subject Classification}: 60 G44.

\

\noindent {\it Keywords}: Stochastic exponential, Girsanov's transformation,
 \\ Local martingale.

\end{abstract}


 {\bf 1. Introduction.} \; Let us introduce a basic probability space $\big( \Omega , \mathcal{F}, P \big)$ and continuous filtration $(\mathcal{F}_t)_{0\leq t \leq \infty}$, which means that every local martingale is continuous. Let $\mathcal{F}_\infty $ be the smallest $\sigma-$Algebra containing all $\mathcal{F}_t $ for $t<\infty$. Let $M = (M_t)_{t \geq 0}$ be a local martingale on the stochastic interval $[[0;T]]$, where $T$ is a stopping time. Denote by $\mathcal E(M)$ the stochastic exponential of a local martingale $M$:

$$
\mathcal E_{t}(M)=\exp\{ M_t-\frac{1}{2}\la M\ra_t \}.
$$
\\
For a given local  martingale $M$, the associated stochastic exponential $\cal {E}(M)$
is a local martingale, but not necessarily a true martingale. To know whether $\cal E(M)$ is a true martingale is important for many applications, e.g., when Girsanov's theorem is  applied to perform a change of measure. 

\
It is well-known that exponential martingales play an essential role in various questions concerning the absolute continuity of probability laws in stochastic processes. A. A. Novikov \cite{6} showed that $\mathcal E (M)$ is a uniformly integrable martingale if $e^{\frac{1}{2}\langle M \rangle_\infty}\in L_1$ and that the constant $\frac{1}{2}$ can not be improved. In 1979 Kazamaki \cite{3} proved that $\sup_\tau Ee^{\frac{1}{2}M_\tau }<\infty $ is sufficient for uniform integrability of $\mathcal E (M)$. Then in 1994 Kazamaki \cite{4} generalized his assertion introducing mixed Novikov-Kazamaki  condition using constant $a\neq 1$ and lower functions (Kazamaki \cite{4}, p.19, Theorem 1.12). In 2013 J. Ruf \cite{7} generalized mixed Novikov-Kazamaki criterion introducing general function of local martingale and its quadratic variation. In \cite{8} and \cite{9} the mixed Novikov-Kazamaki criterion is generalized using predictable process 
$a_s$ instead of the constant $a$. A similar question in the exponential semimartingale framework, in particular, for affine processes, has also attracted attention in Kallsen and Muhle-Kabre \cite{10} and in Kallsen and Shiryaev \cite{11}. In \cite{10} a weak sufficient criterion and in \cite{11} sufficient criterion in terms of cumulant process is given for uniform integrability of $\mathcal{E}(M)$.

\
The necessary and sufficient conditions for the uniform integrability of $\cal E(M)$ were provided in Mayerhofer, Muhle-Kabre and Smirnov \cite{12} by considering the case when the initial martingale $M$ represents one component of a multivariate affine process, and in Blei and Engelbert \cite{13} and Engelbert and Senf \cite{14} for the exponential local martingales associated with a strong Markov continuous local martingale. In \cite{12} deterministic necessary and sufficient conditions is provided in terms of the parameters of the initial martingale $M$. In \cite{13}, the case of a strong Markov continuous local martingale $M$ is studied and the deterministic criterion is expressed in terms of speed measure of $M$. In \cite{14}, the case of a general continuous local martingale $M$ is considered and the condition of uniform integrability of $\cal E(M)$ is given in terms of time-change that turns $M$ into a (possible stopped) Brownian motion. In \cite{17} Yu. M. Kabanov, R. Sh. Liptser, A. N. Shiryaev showed that if the measure $Q$ is locally absolutely continuous w.r.t. the measure $P$, then for absolute continuity of $Q$ w.r.t. $P$ necessary and sufficient is that $Q\{ \la M \ra _T < \infty \} = 1$. For the treatment of the related questions of a local absolute continuity of measures on filtered spaces see also Jacod and Shiryaev \cite{16} and Cheridito, Filipovic and Yor \cite{15}. We establish a necessary and sufficient conditions for the uniform integrability of the stochastic exponential $\cal E(M)$ in terms of the basic measure $P$.   
  
\
\\
In the next section we formulate the main results of this paper ({\bf Theorem 1} and {\bf Theorem 2}) and then we prove them in the third section. In order to prove theorems we need several Lemmas which are given in Appendix.


\

{\bf 2. The main results.} \; In the following theorem we weakened condition $|a_s-1|\geq \varepsilon >0$ imposed in \cite{8} and \cite{9}, which enable us to obtain new type necessary and sufficient condition:

\
\\
{\bf Theorem 1.}
For the uniform integrability of the stochastic exponential $\mathcal{E}(M)$, it is necessary and sufficient, that there exists a predictable, $M-$integrable process $a_s$ such that:
\\
$(i) \;\;\; \sup_{\tau \leq T} E\exp \Big\{\int_0^{\tau } a_s dM_s + \int_0^{\tau } \big( \frac{1}{2}-a_s \big) d \la M \ra _s\Big\} <\infty $
\\
where the $\sup$ is taken over all stopping times $\tau \leq T$;
\\
$(ii) \;\;\; f(\la M \ra _s) \leq (a_s - 1)^2 $ for some function $f\geq 0$ with $\int_0^\infty f(x)dx = \infty $.

\
\\
Notice that because $\mathcal{E}(M)$ is a supermartingale, $E\mathcal{E}_\tau (M)\leq 1$ for any stopping time $\tau \leq T$. So condition $(i)$ of Theorem 1 is automatically satisfied when $a_s\equiv 1$. This means that condition $(i)$ is not sufficient when $a_s$ is quite close to $1$. Accordingly, condition $(ii)$ of Theorem 1 gives us an exact degree of proximity of $a_s$ to $1$.

\
\\
{\bf Remark 1.}  Novikov's \cite{6} condition $Ee^{\frac{1}{2}\langle M \rangle_T}<\infty $ and Kazamaki's \cite{3} criterion $\sup _\tau Ee^{\frac{1}{2} M_\tau } < \infty $ are particular cases of Theorem 1 taking $a_s\equiv 0, \;\; f(x)\equiv 1$ and 
$a_s\equiv \frac{1}{2}, \;\; f(x)\equiv \frac{1}{2}$ respectively. Applying Theorem 1 for $a_s\equiv a\neq 1$ and $f(x)\equiv (1-a)^2$ we obtain the mixed Novikov-Kazamaki's condition: 
$$\sup_{\tau \leq T} Ee^{aM_{\tau } + ( \frac{1}{2}-a ) \la M \ra_{\tau } } < \infty .$$  

\
\\
Let $(B_t)_{t\geq0}$ be a standard Brownian motion. Recall that continuous function $\varphi : R_+ \longrightarrow R_+$ is said to be a {\it lower function} if

$$
P \Big\{ \omega \;\;  : \;\; \exists \; t(\omega), \; \forall \; t>t(\omega)   \Rightarrow  B_t < \varphi(t)   \Big\} = 0.
$$ 
\\
For example, $\varphi(t) = C\sqrt{t}$ and $\varphi(t) = \sqrt{2t\log{\log{t}}}$ are lower functions.

\
\\
In the next theorem we have no restriction on $a_s$, but we have additional term $f( \int^{\tau }_0 1_{\{ |1-a_s|<\varepsilon \}} d\langle M \rangle_s )$ which will be essential when $a_s$ is close to $1$:  

\
\\
{\bf Theorem 2.}\;
For the uniform integrability of the stochastic exponential $\mathcal E(M)$, it is necessary and sufficient, that there exists a predictable process $a_s$, positive constant $\varepsilon $, non-decreasing lower function $\varphi $ and a 
non-decreasing function $f:[0;\infty )\rightarrow (0;\infty )$ with $\lim _{x\rightarrow \infty }f(x)=\infty $ such that:

$$\sup_{\tau \leq T}E\exp \bigg\{ {\int ^\tau_0 a_s dM_s + \int ^\tau_0  \Big( \frac{1}{2} - a_s \Big) d\langle M \rangle_s} - 
\varepsilon \varphi \Big( \int^{\tau }_0 1_{\{ |1-a_s|\geq \varepsilon \}} d\langle M \rangle_s \Big)$$ 

\begin{equation}
+ f\Big( \int^{\tau }_0 1_{\{ |1-a_s|<\varepsilon \}} d\langle M \rangle_s \Big) \bigg\} < \infty 
\end{equation}
where the $\sup$ is taken over all stopping times $\tau \leq T$.

\
\\
{\bf Remark 2.} Novikov's \cite{6} condition $Ee^{\frac{1}{2}\langle M \rangle_T}<\infty $ and Kazamaki's \cite{3} criterion $\sup _\tau Ee^{\frac{1}{2} M_\tau } < \infty $ are particular cases of Theorem 2 taking $a_s\equiv 0$, $\varepsilon = 1$, $\varphi \equiv 0$ and $a_s\equiv \frac{1}{2}$, $\varepsilon = \frac{1}{2}$, $\varphi \equiv 0$ respectively. Applying Theorem 2 for $a_s\equiv a\neq 1$, 
$\varepsilon =|1-a|$ and $\varphi $ lower function we get the mixed Novikov-Kazamaki's condition with non-decreasing lower function:   
$$\sup_{\tau \leq T} E\exp \Big\{aM_{\tau } + \Big( \frac{1}{2}-a \Big) \la M \ra_{\tau } - |1-a|\varphi(\la M \ra_{\tau })\Big\} < \infty .$$ 
Notice that in these cases $f\Big( \int^{\tau }_0 1_{\{ |1-a_s|<\varepsilon \}} d\langle M \rangle_s \Big)=f(0)$.

\
\\
It follows from the proof of necessity of Theorem 2 that Ruf's condition (\cite{7}, Corollary 5) is necessary and sufficient: 

\
\\
{\bf Corollary.}   
For the uniform integrability of the stochastic exponential, it is necessary and sufficient, that there exists a continuous function 
$h:R_+ \rightarrow R$ with $\lim\sup_{x\rightarrow\infty }h(x)=\infty $ such that $\sup _{\tau \leq T} E\mathcal {E}_{\tau }(M)e^{h(\langle M \rangle_{\tau })} < \infty $. 

\
\\
{\bf Remark 3.} It is obvious that Ruf's main result (\cite{7}, {\bf Theorem 1}) is also necessary and sufficient for the uniform integrability of the stochastic exponential.


\

{\bf 3. Proof of the main results.} \;
\\
{\it Proof of Theorem 1. } 
\\
{\it Sufficiency:} According to the condition $(ii)$ from Theorem 1 $\tilde{f}(x)=1+\int_0^x f(t)dt$ \; is a positive and non-decreasing function with 
$\lim_{x\rightarrow \infty }\tilde{f}(x)=\infty $. So using Lemma 1 from appendix there exists absolutely continuous and non-decreasing function 
$\tilde{g}:R_+\rightarrow R_+$ which satisfies conditions $(a)$, $(b)$ and $(c)$ of Lemma 1. Now let us define function $g(x)=\frac{1}{4}\tilde{g}(\frac{x}{2})$. Condition $(b)$ from Lemma 1 implies that $\lim_{x\rightarrow\infty }g(x) = \infty $. Then we will have following inequalities:
$$g\Big(\int_0^{\tau }a_s^2 d\la M \rangle_s \Big)\leq g\Big(2\la M \rangle_{\tau } + 2\int_0^{\tau }(a_s - 1)^2 d\la M \rangle_s \Big)$$     
\begin{equation}
=\frac{1}{4}\tilde{g}\Big(\la M \rangle_{\tau } + \int_0^{\tau }(a_s - 1)^2 d\la M \rangle_s \Big)
\end{equation}
Now according to the conditions $(c)$ and $(a)$ of Lemma 1 from $(2)$ we obtain:
$$g\Big(\int_0^{\tau }a_s^2 d\la M \rangle_s \Big) \leq \frac{1}{4}\tilde{g}\big(\la M \rangle_{\tau } \big) + \frac{1}{4} \int_0^{\tau }(a_s - 1)^2 d\la M \rangle_s + \frac{1}{2}$$
$$\leq \frac{1}{4}\tilde{f}\big(\la M \rangle_{\tau } \big) + \frac{1}{4} \int_0^{\tau }(a_s - 1)^2 d\la M \rangle_s + \frac{1}{2}$$
$$ = \frac{1}{4}\int_0^{\tau } f(\la M \ra_s ) d\la M \ra_s + \frac{1}{4} \int_0^{\tau }(a_s - 1)^2 d\la M \rangle_s + \frac{1}{2}$$
$$\leq \frac{1}{4} \int_0^{\tau }(a_s - 1)^2 d\la M \rangle_s + \frac{1}{4} \int_0^{\tau }(a_s - 1)^2 d\la M \rangle_s + \frac{1}{2} =
 \frac{1}{2} \int_0^{\tau }(a_s - 1)^2 d\la M \rangle_s + \frac{1}{2}.$$
In the last inequality we used condition $(ii)$ from Theorem 1. So as a result we obtained inequality 
$g\Big(\int_0^{\tau }a_s^2 d\la M \rangle_s \Big) \leq \frac{1}{2} \int_0^{\tau }(a_s - 1)^2 d\la M \rangle_s + \frac{1}{2} $, where $g$ is a continuous, non-decreasing function with $\lim_{x\rightarrow\infty }g(x) = \infty $. Now using this last inequality and condition $(i)$ from Theorem 1 we will have:
$$E\mathcal{E}_{\tau }\Big( \int a d M \Big) \exp \Big\{ g \Big( \int_0^{\tau }a_s^2 d\la M \rangle_s \Big) \Big\} $$ 
$$ \leq e^{\frac{1}{2}}E\exp \Big\{\int_0^{\tau }a_s dM_s - \frac{1}{2}\int_0^{\tau } a_s^2 d \la M \ra_s + \frac{1}{2}\int_0^{\tau } (a_s - 1)^2 d \la M \ra_s\Big\}$$

$$=e^{\frac{1}{2}}E\exp \Big\{\int_0^{\tau }a_s dM_s + \int_0^{\tau } \Big( \frac{1}{2} - a_s \Big) d \la M \ra_s\Big\} < \infty .$$
According to Ruf's condition (\cite{7}, Corollary 5) this implies that $\mathcal{E}(\int a d M)$ is a uniformly integrable martingale. So we have the equality $E\mathcal{E}_T(\int a d M) = 1$.  
\\
Now define the new probability measure $d\tilde{P}=\mathcal{E}_T(\int adM)dP$ and local martingale $N_t=\int^t_0(1-a_s)dM_s$. According to Girsanov's theorem 
$$\tilde{N}_t = N_t - \langle N;\int adM \rangle_t $$
$$= \int^t_0(1-a_s)dM_s - \int^t_0 a_s(1-a_s)d\langle M \rangle_s $$
is a $\tilde{P}-$local martingale. Let us show that for $\tilde{N}$ Novikov's condition is satisfied:
$$E^{\tilde{P}}\exp \Big\{\frac{1}{2}\langle \tilde{N} \rangle_T \Big\} = E \exp \Big\{ \int^T_0 a_sdM_s -\frac{1}{2}\int^T_0 a_s^2 d\langle M \rangle_s +\frac{1}{2}\int^T_0 (1-a_s)^2 d\langle M \rangle_s \Big\}$$

$$= E\exp \Big\{ \int^T_0 a_sdM_s + \int^T_0 \Big(\frac{1}{2} - a_s \Big) d\langle M \rangle_s \Big\} < \infty $$
by condition $(i)$ from Theorem 1.
This means that $E^{\tilde{P}}\mathcal{E}_T(\tilde{N})=1$. Finally we get 
$$1=E^{\tilde{P}}\mathcal{E}_T(\tilde{N})=E\exp \Big\{ \int^T_0 a_sdM_s -\frac{1}{2}\int^T_0 a_s^2 d\langle M \rangle_s \Big\} $$
$$\times \exp \Big\{ \int^T_0 (1-a_s)dM_s - \int^T_0 a_s(1-a_s) d\langle M \rangle_s - \frac{1}{2}\int^T_0 (a_s-1)^2 d\langle M \rangle_s \Big\} $$
$$= E\exp\Big\{ M_T - \frac{1}{2}\int^T_0 a_s^2 d\langle M \rangle_s - \int^T_0 a_s d\langle M \rangle_s + \int^T_0 a_s^2 d\langle M \rangle_s $$
$$ - \frac{1}{2}\int^T_0 a_s^2 d\langle M \rangle_s + \int^T_0 a_s d\langle M \rangle_s -\frac{1}{2}\langle M \rangle_T \Big\} =E\mathcal{E}_T(M). $$
Thus $E\mathcal{E}_T(M)=1$, which implies that $\mathcal{E}(M)$ is a uniformly integrable martingale.

\
\\
{\it Necessity:}\;\;\; Now we know that $\mathcal{E}(M)$ is a uniformly integrable martingale which implies that $E\mathcal{E}_T(M)=1$. So we can define new probability measure $dQ=\mathcal{E}_T(M)dP$ and a $Q-$local martingale 
$\tilde{M} = M - \la M \ra$. It follows from \cite{17} that $Q(\langle \tilde{M} \rangle_T < \infty )=1$, so according to Lemma 2 there exists absolutely continuous and non-decreasing function $h$ with $\lim_{x\rightarrow \infty }h(x)=\infty $ such that $E^Q e^{h(\la \tilde{M} \ra_T)} < \infty $. Because $h$ is non-decreasing and absolutely continuous, there exists non-negative function $f$ such that $h(x) = h(0) + \int_0^x f(t)dt$ and therefore $h(\la \tilde{M} \ra_\tau) = h(0) + \int_0^\tau f(\la \tilde{M} \ra_s)d\la \tilde{M} \ra_s$. So if we define $a_s$ such that $2(a_s - 1)^2 = f(\la \tilde{M} \ra_s)$ then the condition $(ii)$ of Theorem 1 will be satisfied. Now Let us check condition $(i)$:
 
$$E\exp \Big\{ \int^\tau _0 a_sdM_s + \int^\tau _0 \Big(\frac{1}{2} - a_s \Big) d\langle M \rangle_s \Big\} 
= E^Q \exp \Big\{ \int^\tau _0 (a_s-1)d\tilde{M}_s \Big\}$$
$$
= E^Q \exp \Big\{ \int^\tau _0 (a_s-1)d\tilde{M}_s - \int^\tau _0 (a_s - 1)^2 d \la \tilde{M} \ra_s \Big\}\times \exp \Big\{ \int^\tau _0 (a_s - 1)^2 d \la \tilde{M} \ra_s \Big\} 
$$

$$\leq \sqrt{E^Q \mathcal{E}_\tau \Big( 2\int (a-1) d \tilde{M} \Big)} \times \sqrt{E^Q \exp \Big\{ 2\int^T_0 (a_s - 1)^2 d \la \tilde{M} \ra_s \Big\}}$$

$$
\leq \sqrt{E^Q \exp \Big\{ 2\int^T_0 (a_s - 1)^2 d \la \tilde{M} \ra_s \Big\}} = \sqrt{E^Q \exp \Big\{ h(\la \tilde{M} \ra_T) - h(0) \Big\}} < \infty .
$$
Proof of Theorem 1 is completed.
\qed

\
\\
{\it Proof of Theorem 2.} 
\\
{\it Sufficiency:}
Let us first show that $E\mathcal{E}_T(\int adM)=1$. It follows from Ruf's condition (\cite{7}, Corollary 5) that for this it is sufficient to show that 
$$
\sup _{\tau \leq T}E\mathcal{E}_{\tau }\Big(\int adM \Big)\exp \Big\{ h\Big( \int^{\tau }_0 a_s^2 d\langle M \rangle_s \Big) \Big\} < \infty
$$
for some continuous function $h$ with $\limsup _{x\rightarrow\infty }h(x)=\infty $. According to the conditions of Theorem 2 function $f$ from $(1)$ is non-decreasing, so using Lemma 1 there exists positive, non-decreasing and absolutely continuous function $g$ which satisfies conditions $(a)$, $(b)$ and $(c)$ of Lemma 1. Let us define function $h$ as $h(x)=\delta g\big( \frac{1}{(1+\varepsilon )^2} x\big)$ where $0<\delta <1$ is a constant sufficiently close to $0$. It is well known that $\limsup _{x\rightarrow\infty }\frac{\varphi (x)}{x}=0$ for any non-decreasing lower function $\varphi $, which implies inequality $\varphi(x)\leq \delta x + C_\delta $ for any $\delta >0$. Now using the definition of $h$ and inequality $\varphi(x)\leq \delta x + C_\delta $ we get:

$$ h\Big(\int^{\tau }_0 a_s^2 d\langle M \rangle_s\Big) + \varepsilon \varphi \Big( \int^{\tau }_0 1_{\{ |1-a_s|\geq \varepsilon \}} d\langle M \rangle_s \Big) $$
$$\leq \delta g\Big(\frac{1}{(1+\varepsilon )^2}\int^{\tau }_0 a_s^2 d\langle M \rangle_s\Big) + \varepsilon \delta \int^{\tau }_0 1_{\{ |1-a_s|\geq \varepsilon \}} d\langle M \rangle_s + \varepsilon C_\delta $$
$$ = \delta g\bigg(\frac{1}{(1+\varepsilon )^2}\int^{\tau }_0 a_s^2 1_{\{|a_s-1|<\varepsilon \}} d\langle M \rangle_s + \frac{1}{(1+\varepsilon )^2}\int^{\tau }_0 a_s^2 1_{\{|a_s-1|\geq \varepsilon \}} d\langle M \rangle_s\bigg) $$
\begin{equation}
+ \varepsilon \delta \int^{\tau }_0 1_{\{ |1-a_s|\geq \varepsilon \}} d\langle M \rangle_s + \varepsilon C_\delta 
\end{equation}
Now using the inequality $g(x+y) \leq g(x) + y + 2$  we will have from $(3)$:
$$ h\Big(\int^{\tau }_0 a_s^2 d\langle M \rangle_s\Big) + \varepsilon \varphi \Big( \int^{\tau }_0 1_{\{ |1-a_s|\geq \varepsilon \}} d\langle M \rangle_s \Big) $$
$$\leq \delta g\bigg( \frac{1}{(1+\varepsilon )^2}\int^{\tau }_0 a_s^2 1_{\{|a_s-1|<\varepsilon \}} d\langle M \rangle_s \bigg) + \frac{\delta }{(1+\varepsilon )^2}\int^{\tau }_0 a_s^2 1_{\{|a_s-1|\geq \varepsilon \}} d\langle M \rangle_s  $$

\begin{equation}
+ 2\delta + \varepsilon \delta \int^{\tau }_0 1_{\{ |1-a_s|\geq \varepsilon \}} d\langle M \rangle_s + \varepsilon C_\delta 
\end{equation}
It is clear that $ \frac{1}{(1+\varepsilon )^2}\int^{\tau }_0 a_s^2 1_{\{|a_s-1|<\varepsilon \}} d\langle M \rangle_s \leq \int^{\tau }_0 1_{\{|a_s-1|<\varepsilon \}} d\langle M \rangle_s$ and 
$$\frac{\delta }{(1+\varepsilon )^2}\int^{\tau }_0 a_s^2 1_{\{|a_s-1|\geq \varepsilon \}} d\langle M \rangle_s + \varepsilon \delta \int^{\tau }_0 1_{\{ |1-a_s|\geq \varepsilon \}} d\langle M \rangle_s $$
$$\leq \int^{\tau }_0 (\delta a_s^2 + \delta \varepsilon )1_{\{ |1-a_s|\geq \varepsilon \}} d\langle M \rangle_s .$$
So if we use inequality $g(x)\leq f(x)$ and non-decreasing property of $f$ we obtain from $(4)$:
$$ h\Big(\int^{\tau }_0 a_s^2 d\langle M \rangle_s\Big) + \varepsilon \varphi \Big( \int^{\tau }_0 1_{\{ |1-a_s|\geq \varepsilon \}} d\langle M \rangle_s \Big) $$

$$
\leq f\Big( \int^{\tau }_0 1_{\{|a_s-1|<\varepsilon \}} d\langle M \rangle_s \Big) + \int^{\tau }_0 (\delta a_s^2 + \delta \varepsilon)1_{\{|a_s-1|\geq \varepsilon \}} d\langle M \rangle_s + 2\delta + \varepsilon C_\delta
$$
\begin{equation}
\leq f\Big( \int^{\tau }_0 1_{\{|a_s-1|<\varepsilon \}} d\langle M \rangle_s \Big) + \frac{1}{2} \int^{\tau }_0 (a_s-1)^2 d\langle M \rangle_s + 2\delta + \varepsilon C_\delta .
\end{equation}
In the last inequality we used Lemma 4 to obtain estimation $(\delta a_s^2 + \delta \varepsilon)1_{\{|a_s-1|\geq \varepsilon \}} \leq \frac{1}{2}(a_s -1)^2$ for 
$\delta >0$ sufficiently close to $0$.
So in $(5)$ we got the inequality
$$h\Big(\int^{\tau }_0 a_s^2 d\langle M \rangle_s\Big) + \varepsilon \varphi \Big( \int^{\tau }_0 1_{\{ |1-a_s|\geq \varepsilon \}} d\langle M \rangle_s \Big)$$
$$\leq f\Big( \int^{\tau }_0 1_{\{|a_s-1|<\varepsilon \}} d\langle M \rangle_s \Big) + \frac{1}{2} \int^{\tau }_0 (a_s-1)^2 d\langle M \rangle_s + 2\delta +C_\delta $$
which is equivalent to the following
$$-\frac{1}{2}\int^{\tau }_0 a_s^2 d\langle M \rangle_s + h\Big(\int^{\tau }_0 a_s^2 d\langle M \rangle_s\Big) \leq \int^{\tau }_0 \Big( \frac{1}{2} - a_s \Big)d\langle M \rangle_s$$
\begin{equation}
-\varepsilon \varphi \Big( \int^{\tau }_0 1_{\{ |1-a_s|\geq \varepsilon \}} d\langle M \rangle_s \Big)+ f\Big( \int^{\tau }_0 1_{\{|a_s-1|<\varepsilon \}} d\langle M \rangle_s \Big) + 2 + \varepsilon C_\delta .
\end{equation}
Now from $(6)$ we obtain 
$$ \sup _{\tau \leq T}E\mathcal{E}_{\tau }\Big(\int adM \Big)\exp \Big\{ h\Big( \int^{\tau }_0 a_s^2 d\langle M \rangle_s \Big) \Big\} $$
$$\leq e^{2 + \varepsilon C_\delta }\sup_{\tau \leq T}E\exp \bigg\{ {\int ^\tau_0 a_s dM_s + \int ^\tau_0  \Big( \frac{1}{2} - a_s \Big) d\langle M \rangle_s} $$
$$- \varepsilon \varphi \Big( \int^{\tau }_0 1_{\{ |1-a_s|\geq \varepsilon \}} d\langle M \rangle_s \Big)+ f\Big( \int^{\tau }_0 1_{\{ |1-a_s|<\varepsilon \}} d\langle M \rangle_s \Big) \bigg\} < \infty .$$
According to Ruf's condition (\cite{7}, Corollary 5) this means that $\mathcal{E}(\int adM)$ is a uniformly integrable martingale which implies that
$E\mathcal{E}_T(\int adM)=1 $.  
\
\\
Now define the new probability measure $d\tilde{P}=\mathcal{E}_T(\int adM)dP$ and local martingale $N_t=\int^t_0(1-a_s)dM_s$. According to Girsanov's theorem 
$$\tilde{N}_t = N_t - \langle N;\int adM \rangle_t $$
$$= \int^t_0(1-a_s)dM_s - \int^t_0 a_s(1-a_s)d\langle M \rangle_s $$
is a $\tilde{P}-$local martingale. Define function $\psi(x)=\varepsilon \varphi (\frac{x}{\varepsilon^2})$ which is lower function according to Lemma 3. Let us show that for $\tilde{N}$ Novikov's condition with lower function is satisfied:
$$E^{\tilde{P}}\exp \Big\{\frac{1}{2}\langle \tilde{N} \rangle_\tau - \psi (\langle \tilde{N} \rangle_\tau )\Big\} = E \exp \Big\{ \int^\tau_0 a_sdM_s -\frac{1}{2}\int^\tau_0 a_s^2 d\langle M \rangle_s $$
$$+\frac{1}{2}\int^\tau_0 (1-a_s)^2 d\langle M \rangle_s - \varepsilon \varphi \Big(\frac{1}{\varepsilon^2}\int^\tau_0 (1-a_s)^2 d\langle M \rangle_s\Big) \Big\}$$
$$
= E\exp \Big\{ \int^\tau_0 a_sdM_s + \int^\tau_0 \Big(\frac{1}{2} - a_s \Big) d\langle M \rangle_s - \varepsilon \varphi \Big(\frac{1}{\varepsilon^2}\int^\tau_0 (1-a_s)^2 d\langle M \rangle_s\Big)\Big\}
$$
$$\leq E\exp \Big\{ \int^\tau_0 a_sdM_s + \int^\tau_0 \Big(\frac{1}{2} - a_s \Big) d\langle M \rangle_s $$
\begin{equation}
- \varepsilon \varphi \Big(\int ^\tau _0 1_{\{ |1-a_s|\geq \varepsilon \}} d\langle M \rangle_s \Big)\Big\}< \infty .
\end{equation}
Here we used the inequality 
$$\varepsilon \varphi \Big( \int^{\tau }_0 1_{\{ |1-a_s|\geq \varepsilon \}} d\langle M \rangle_s \Big) \leq \varepsilon \varphi \Big(\frac{1}{\varepsilon^2}\int^\tau_0 (1-a_s)^2 d\langle M \rangle_s\Big) $$
which follows from non-decreasing property of lower function $\varphi $. 
\\
$(7)$ implies that $E^{\tilde{P}}\mathcal{E}_T(\tilde{N})=1$ which is equivalent to the $E\mathcal{E}_T(M)=1$: 
$$1=E^{\tilde{P}}\mathcal{E}_T(\tilde{N})=E\exp \Big\{ \int^T_0 a_sdM_s -\frac{1}{2}\int^T_0 a_s^2 d\langle M \rangle_s \Big\} $$
$$\times \exp \Big\{ \int^T_0 (1-a_s)dM_s - \int^T_0 a_s(1-a_s) d\langle M \rangle_s - \frac{1}{2}\int^T_0 (a_s-1)^2 d\langle M \rangle_s \Big\} $$
$$= E\exp\Big\{ M_T - \frac{1}{2}\int^T_0 a_s^2 d\langle M \rangle_s - \int^T_0 a_s d\langle M \rangle_s + \int^T_0 a_s^2 d\langle M \rangle_s $$
$$ - \frac{1}{2}\int^T_0 a_s^2 d\langle M \rangle_s + \int^T_0 a_s d\langle M \rangle_s -\frac{1}{2}\langle M \rangle_T \Big\} =E\mathcal{E}_T(M). $$
Thus $E\mathcal{E}_T(M)=1$, which implies that $\mathcal{E}(M)$ is a uniformly integrable martingale.

\
\\
{\it Necessity:} \;\;\; 
Now we know that $\mathcal{E}(M)$ is a uniformly integrable martingale. So, we can define the new probability measure $dQ=\mathcal{E}_T(M)dP$. It follows from \cite{17} that $Q(\la M \ra _T < \infty ) = 1$, so we can apply Lemma 2 to find non-decreasing and absolutely continuous function $f$ with 
$\lim _{x \rightarrow \infty }f(x)=\infty $ such that $E^Qe^{f(\la M \ra _T)} < \infty $. Now if we insert $a_s\equiv 1$, $\varphi \equiv 0$, $\varepsilon >0$ and $f$ in $(1)$ we obtain:
$$\sup_{\tau \leq T}E\exp \bigg\{ {\int ^\tau_0 a_s dM_s + \int ^\tau_0  \Big( \frac{1}{2} - a_s \Big) d\langle M \rangle_s} - 
\varepsilon \varphi \Big( \int^{\tau }_0 1_{\{ |1-a_s|\geq \varepsilon \}} d\langle M \rangle_s \Big)$$ 

$$ + f\Big( \int^{\tau }_0 1_{\{ |1-a_s|<\varepsilon \}} d\langle M \rangle_s \Big) \bigg\} = \sup _{\tau \leq T} E\exp \bigg\{ M_\tau -\frac{1}{2}\la M\ra_\tau - \varepsilon \varphi (0)+ f( \langle M \rangle_\tau ) \bigg\}$$
$$ \leq \sup _{\tau \leq T} E\mathcal{E}_\tau (M)e^{f(\la M \ra _\tau)} = \sup _{\tau \leq T} E^Q e^{f(\la M \ra _\tau)} \leq E^Q e^{f(\la M \ra _T)} < \infty .$$   
Proof of Theorem 2 is completed.

\qed


\

{\bf 4. Appendix.} \;
\\
{\bf Lemma 1.} \; Let $f:[0;+\infty ) \rightarrow (0; +\infty )$ be a non-decreasing function with $\lim _{x\rightarrow \infty }f(x) = \infty $. Then 
a non-decreasing, absolutely continuous function $g:[0;+\infty ) \rightarrow [0; +\infty )$ exists which satisfies the following conditions:
\\
$(a) \;\;\; g(x)\leq f(x)$;
\\
$(b) \;\;\; \lim _{x\rightarrow \infty }g(x) = \infty $;
\\
$(c) \;\;\; g(x+y)\leq g(x)+y+2$.
\
\\
{\it proof.} \;        
Define function $F$: $F(x) = \sum^{\infty }_{k=1}f(k-1)1_{[k-1; k[}(x)$. It is obvious that $F(x)\leq f(x)$. Let us denote by 
$\bigtriangleup F(k)=F(k)-F(k-1)$ jumps of $F$. Because $f$ is non-decreasing, the jumps of $F$ will be non negative. Now define non-decreasing sequence  
$(g_k)_{k\geq 1}$ with recurrence relationship:    
$$g_0=0; \;\; g_1=1\wedge f(0); \;\; g_2=g_1+1\wedge \bigtriangleup F(2); \;\;\; g_k=g_{k-1}+1\wedge \bigtriangleup F(k) \;\;\;\;\; k\geq 2.$$
As a result we have points $(k; g_k) \;\;\; k\geq 0$. Define function $g$ by connecting points $(k; g_k)$ with straight lines. It follows from the definition that $g$ is absolutely continuous, non-decreasing, $g(x)\leq F(x) \leq f(x)$ and $\lim _{x\rightarrow \infty }g(x) = \infty $. It remains to show that $g(x+y)\leq g(x)+y+2$. Let us take
$x\in [k-1;k]$ and $y\in [n-1;n]$. It is clear that $x+y \leq k+n$. Using the definition of function $g$ we obtain:
$$g(x+y)-g(x)\leq g(k+n)-g(k-1) = \sum _{i=1}^{k+n}1\wedge \bigtriangleup F(i) - \sum _{i=1}^{k-1}1\wedge \bigtriangleup F(i) = $$
$$ = \sum _{i=k}^{k+n}1\wedge \bigtriangleup F(i) \leq n+1 \leq y+2.$$  
So, by arbitrariness of $k$ and $n$, $g(x+y)\leq g(x)+y+2$.
\qed 

\
\\
{\bf Lemma 2.} \; For any random variable $\xi $ such that $P(0\leq \xi < \infty )=1$ there exists a positive, non-decreasing and continuous function $g$ with 
$\lim _{x\rightarrow \infty }g(x) = \infty $, such that $Eg(\xi ) < \infty $.
\
\\
{\it Proof.} \; Let $F_{\xi }(x) = P(\xi \leq x)$ be the probability distribution function of $\xi $. Let us take $f(x)=\frac{1}{\sqrt{1 - F_{\xi }(x-)}}$. Then we will have:
$$ Ef(\xi ) = \int^{\infty }_0 \frac{1}{\sqrt{1 - F_{\xi }(x-)}} dF_{\xi }(x) = \bigg[ -2\sqrt{1-F_{\xi }(x)} \bigg]^{\infty }_0  $$
$$  - \sum_{0<x<\infty }\Big[ -2\sqrt{1-F_{\xi }(x)} + 2\sqrt{1-F_{\xi }(x-)}-\frac{\Delta F_{\xi }(x)}{\sqrt{1-F_{\xi }(x-)}} \Big]  $$
$$ \leq \bigg[ -2\sqrt{1-F_{\xi }(x)} \bigg]^{\infty }_0 = 2\sqrt{1-F_{\xi }(0)} < \infty .$$
Here we have used inequality $-2\sqrt{1-F_{\xi }(x)} + 2\sqrt{1-F_{\xi }(x-)}-\frac{\Delta F_{\xi }(x)}{\sqrt{1-F_{\xi }(x-)}} \geq 0$ which follows from the convexity of the function $y \rightarrow -2\sqrt{1-y}$. Now if we use Lemma 1 we can find absolutely continuous, positive and non-decreasing function $g$ with 
$\lim _{x\rightarrow \infty }g(x) = \infty $, such that $g(x)\leq f(x)$. Inequalities $g(x)\leq f(x)$ and $Ef(\xi ) < \infty $ implies that $Eg(\xi ) < \infty $.  
\qed

\
\\
{\bf Lemma 3.}\;\;\;
Let $\varphi $ be a lower function. Then the function $\psi(x)=\varepsilon\varphi(\frac{x}{\varepsilon^2})$ also will be a lower function. 
\
\\
{\it Proof. }  \;\;\; It is well known that if $B_t$ is a Brownian motion, then $W_t = \varepsilon B_{t/\varepsilon^2}$ will be Brownian motion too. With this if we take $s=t/\varepsilon^2$ then we will have:
$$
P \Big\{ \omega \;\;  : \;\; \exists \; t(\omega), \; \forall \; t>t(\omega)   \Rightarrow  W_t < \psi(t)   \Big\} =
$$
$$
= P \Big\{ \omega \;\;  : \;\; \exists \; s(\omega), \; \forall \; s>s(\omega)   \Rightarrow  B_s < \varphi(s)   \Big\} = 0.
$$ 
\qed

\
\\
{\bf Lemma 4.}\;\;\;
For any $\varepsilon >0$ there exists $\delta > 0$ sufficiently close to $0$ such that inequality $\delta x^2 + \delta \varepsilon \leq \frac{1}{2}(x-1)^2$ holds true for any $x \notin (1-\varepsilon ; 1+\varepsilon )$.
\
\\
{\it proof. } \; It is clear that we can take $\delta >0$ sufficiently close to $0$ such that 
$$1-\varepsilon < \frac{1-\sqrt{2\delta (1+\varepsilon -2\delta \varepsilon)}}{1-2\delta } < \frac{1+\sqrt{2\delta (1+\varepsilon -2\delta \varepsilon)}}{1-2\delta } < 1+\varepsilon .$$
It is easy to check that for such $\delta $ condition $|x-1|\geq \varepsilon $ implies inequality 
$$(1-2\delta )x^2 - 2x + 1-2\delta \varepsilon \geq 0$$
which is equivalent to the following $\delta x^2 + \delta \varepsilon \leq \frac{1}{2}(x-1)^2$.
\qed

\end{document}